\documentclass[reqno,a4paper]{elsart}

\usepackage{latexsym}  
\usepackage{exscale}
\usepackage{natbib}                    
\usepackage{amsmath,amsfonts,amssymb}

\newtheorem{theorem}{Theorem}
\newtheorem{lemma}[theorem]{Lemma}

\newcommand{\defi}{\hbox{\small{$\;\stackrel{\text{def}}{=}\;$}}}
\renewcommand{\liminf}{\operatornamewithlimits{\varliminf}}
\renewcommand{\limsup}{\operatornamewithlimits{\varlimsup}}

\newtheorem{teora}{Theorem}
 
\newtheorem{corola}[teora]{Corollary}

\begin{document}
\begin{frontmatter}

\title{Large Deviations and Phase Transition for Random Walks in   
Random Nonnegative Potentials}

\author{Markus Flury}

\address{Universit\"at Z\"urich, Institut f\"ur Mathematik,
Winterthurerstr.~190, CH-8057 Z\"urich, Switzerland}
\ead{mflury@amath.unizh.ch}

\begin{abstract}
We establish 
large deviation
principles  and   phase transition results  for both
quenched and annealed settings of  nearest-neighbor random walks with
constant drift in  random nonnegative potentials on
$\mathbb Z^d$. 
We complement the analysis of \cite{Zer}, where a
shape theorem on the Lyapunov functions and a large deviation
principle in absence of the drift are achieved for the  quenched
setting.
\end{abstract}


\end{frontmatter}

\section{Introduction}\label{Introduction}

Let $\mathcal S=(S(n))_{n\in\mathbb{N}_0}$ be a symmetric nearest-neighbor
random walk on
$\mathbb Z^d$  starting at the origin, and denote by $P$, respectively  $E$,
the associated probability measure, respectively expectation. The aim of this
article is a probabilistic description of  the long-time behavior of the random
walk, endowed with a drift and evolving in  a random environment given by a
random potential on the lattice. This description will be done for concrete
realizations of the environment, the \emph{quenched} setting, as
well as for the averaged environment, the so-called
\emph{annealed} setting. For details, we make the following assumptions\,:
\begin{enumerate}
\item[Qu)]
$\mathbb V=(V_x)_{x\in\mathbb Z^d}$ is a family of independent, identically and
not trivially distributed
random variables in $L^d(\Omega, \mathcal F,\mathbb P)$, which is independent of 
the random walk itself and satisfies  
$\mathrm{ess\,inf}\,V_x=0$.

\item[An)]
$\varphi:[0,\infty)\to[0,\infty)$ is a non constant, non decreasing and concave 
 function with $\varphi(0)=0$ and
$\lim_{t\to\infty}\varphi(t)/t=0$.
\end{enumerate}

For 
$\omega\in\Omega$, 
$n\in\mathbb N$ and
$h\in\mathbb R^d$, the \emph{ quenched path measure
$Q^h_{n,\omega}$} for the random walk $\mathcal S$ with constant
drift~$h$ under the path  potential 
$$\Psi(n,\omega)\defi\sum_{1\leq m\leq n}V_{S(m)}(\omega)$$
is defined by means of the density function
$$\frac{d Q^h_{n,\omega}}{d P}
\defi\frac{1}{Z^h_{n,\omega}}\,
\exp\big(h\cdot S(n)-\Psi(n,\omega)\big)\,,$$
where $Z^h_{n,\omega}$
denotes the corresponding (quenched) normalization.
Notice that  $Q^h_{n,\cdot}$ is a random probability measure, the
randomness coming from the random potential $\Psi(n,\cdot)$.
For $x\in\mathbb Z^d$, let now
$$l_{x}(n)\defi \sharp \left\{1\leq m\leq n:\; S(m) =x\right\}$$ denote the
number of the random walk's
visits to the site $x$ up to time $n$.
The \emph{annealed
path measure
$Q^h_n$} for the random walk $\mathcal S$  with constant drift~$h$
under the path potential
$$
\Phi(n)\defi\sum_{x\in\mathbb Z^d}
\varphi(l_x(n))
$$
is defined by means of the density function
$$
\frac{d Q^h_n}{d P}\defi\frac{1}{Z^h_n}\,
\exp\big(h\cdot S(n)-\Phi(n)\big)\,,
$$
where $Z^h_n$ is the corresponding (annealed) normalization
constant.

The model we come to introduce is a discrete-setting model for a particle moving
in a random media.
In the quenched setting, the walker jumps from site to site, thereby trying to
avoid those regions where the potential takes on high values. 
The drift however implies a restriction in the search of such an ``optimal
strategy'' by imposing a particular direction to the walk. 

We shall point out that in the definition of the
annealed path measures we are  making  a slight abuse of standard
terminology. 
To clarify this aspect,
consider 
\begin{equation*}
\varphi_\mathbb V(t) \defi-\log\mathbb E\, \exp(-t V_x)\,,\quad
t\in[0,\infty)\,,
\end{equation*} 
for a given potential
$\mathbb V$.
By H\"older inequality, 
dominated convergence and the assumption
$\mathrm{ess\,inf}\,V_x=0$, it is easy to see that  $\varphi_\mathbb
V$ fulfills the requirements  An).  
Let $Q^h_{\mathbb V,n}$ denote the annealed path measure corresponding to
$\varphi_\mathbb V$.
The quenched  potential obviously can
be rewritten as
$$
\Psi(n,\omega)=\sum_{x\in\mathbb Z^d} l_x(n) V_x(\omega)\,.$$
By  the  independence assumption on the potential, it now is easily seen that
$$\frac{d Q^h_{\mathbb V,n}}{d P}=
\frac{1}{\mathbb E\, Z^h_{n,\cdot}}\,
\mathbb E\left[\exp\big(h\cdot S(n)-\Psi(n,\cdot)\big)\right]
$$
for any drift~$h$ and all  $n\in\mathbb N$, which is the 
``classical'' annealed path measure.
Our results   cover this standard case, but do not rely on the particular form
of~$\varphi_\mathbb V$ in the above definition.

An interesting example of such a potential $\varphi_\mathbb V$  is
considered at the so-called
\emph{hard obstacle} or \emph{trap model}. There, one assumes
$V_x:\Omega\to\{0,\infty\}$ with positive probability for both values.  The
name of the model comes from the fact that
$$
Z^h_{n,\omega}=E\,\big[\exp\left(h\cdot S(n)\right);\,\text{$V_{S(m)}(\omega)=0$
for
$1< m\leq n$}\big]\,,$$ 
which describes the probability for the drifted random walk not to step into
one of the  ``traps'' $\{x\in\mathbb Z^d:\, V_x(\omega)=\infty\}$ up to time
$n\in\mathbb N$\,. Such a potential is not in $L^1(\Omega,\mathcal F,\mathbb P)$
and  consequently does not satisfy  assumption Qu).   Yet, the function
$\varphi_\mathbb V(t)$, associated to the classical annealed terms, does 
fulfill the required properties An), and  satisfies
$$\varphi_\mathbb V(t)=
\left\{
\begin{array}{cl}
-\log\mathbb P\left[\,V_x=0\,\right]& \text{if $t>0$}\,,\\
0& \text{if $t=0$}\,.
\end{array}\right.$$
Again,  the expected probability $\mathbb E\,
Z^h_{n,\cdot}$ of not stepping into a trap equals the annealed normalization
constant $Z^h_{\mathbb V,n}$ corresponding to $\varphi_\mathbb V$. 
We thus have
\begin{align*}
\mathbb E\, Z^h_{n,\cdot}
 =E\big[\exp\big(h\cdot S(n)-\gamma\,\sharp\{S(m):\, 1\leq m\leq
n\}\,\big)\big]
\end{align*}
with $\gamma=-\log\mathbb P[\,V_x=0\,]$ and any $n\in\mathbb N$.

We come back to the general setup of a random walk in a random potential.
A similar model in a continuous setting, namely  Brownian motion in a
Poissonian potential, was first studied by A.S.~Sznitman. By means of the
powerful method of \emph{enlargement of obstacles}, Sznitman  established a
precise picture in both  quenched and annealed settings. He achieved results
such as a shape theorem, large deviation principles (LDP's) and an accurate
description of the transition between small and large drift.  We refer the
reader to Chapter~5 of \cite{Szn} for a complete review of these results. 
In the discrete setting, an ample study of the random walk under the
influence of the  quenched potential was made by \cite{Zer}.
His results, however, are limited to the case where no drift is present.

The aim of the present work is to add the missing pieces to Zerner's
analysis, recovering the larger picture for the random walk with drift in both
quenched and annealed settings.
The organization of the article is as follows\,:
In Section~\ref{Main_results}, we state the main results. 
In Section~\ref{Shape_theorem}, we  follow Zerner's analysis and  prove a shape
theorem for the directed random walk.
Section~\ref{Large_deviation_principles}
is devoted to the proof of the LDP's. 
In Section~\ref{dual_norms_phase_transitions}, we closely follow Sznitman's
path to analyze  phase transitions in the long-time behavior of the random walk,
related to the size of the drift.

\section{Main results}\label{Main_results}

The essential quantities in our  study of the large  time asymptotics of
the random walk are the so-called \emph{Lyapunov functions} on $\mathbb R^d$.
Let  
$$H(x)\defi\inf\left\{n\in\mathbb N_0:\, S(n)=x\right\}$$
denote the  time of the random walk's first visit to 
the lattice site $x\in\mathbb Z^d$.
For $\lambda\geq 0$, $\omega\in\Omega$ and $x\in\mathbb Z^d$, we define the two-point functions
\begin{align*}
a_\lambda(x,\omega)&\defi-\log E\left[\exp\big(-\lambda
H(x)-\Psi(H(x),\omega)\big);\, H(x)<\infty\right]\,,\\
b_\lambda(x)&\defi-\log E\left[\exp\big(-\lambda
H(x)-\Phi(H(x))\big);\, H(x)<\infty\right]\notag\,.
\end{align*}

Our first result introduces the Lyapunov functions $\alpha_\lambda$ and
$\beta_\lambda$, and sets them in relation to the asymptotic behavior of
$a_\lambda$ and $b_\lambda$.

\begin{teora}[Shape theorem]\label{shape_theorem}\noindent
\begin{enumerate}
\item[a)] {\citep{Zer}}
There is a family
$(\alpha_\lambda)_{\lambda\geq 0}$ of norms on
$\mathbb R^d$  such that for any $\lambda\geq 0$ and all sequences
$(x_k)_{k\in\mathbb N}$ on $\mathbb Z^d$ with
$\|x_k\|_1\to\infty$ as $k\to\infty$, we have 
\begin{equation}\label{quenched_shape_limit}
\lim_{k\to\infty}\,\frac{a_\lambda(x_k,\omega)}{\alpha_\lambda(x_k)}=1
\end{equation}
on a set $\Omega_\lambda$ of full $\mathbb P$-measure and in
$L^1(\Omega,\mathcal F,\mathbb P)$.  Moreover, $\alpha_{\lambda}(x)$ is
continuous in
$(\lambda,x)\in[0,\infty)\times \mathbb R^d$,  concave increasing in
$\lambda\in[0,\infty)$, and satisfies
\begin{equation}\label{bounds_alpha}
{\|x\|_1}\left(\lambda-\log\mathbb E\,\exp(-V_x)\right)\leq{\alpha_\lambda(x)}
\leq{\|x\|_1}\left(\lambda+\log(2 d)+\mathbb EV_x\right)\,.
\end{equation}\pagebreak

\item[b)]
There is a family $(\beta_\lambda)_{\lambda\geq 0}$ of norms on $\mathbb R^d$
such that for any $\lambda\geq 0$ and all sequences
$(x_k)_{k\in\mathbb N}$ on $\mathbb Z^d$ with
$\|x_k\|_1\to\infty$ as $k\to\infty$,
we have
\begin{equation}\label{annealed_shape_limit}
\lim_{k\to\infty}\,\frac{b_\lambda(x_k)}{\beta_\lambda(x_k)}=1\,.
\end{equation}
Moreover, $\beta_{\lambda}(x)$ is continuous in
$(\lambda,x)\in[0,\infty)\times \mathbb R^d$,  concave increasing in
$\lambda\in[0,\infty)$, and satisfies
\begin{equation*}
\quad
{\|x\|_1}\left(\lambda+\varphi(1)\right)\leq{\beta_\lambda(x)}
\leq {\|x\|_1}\left(\lambda+\log(2 d)+\varphi(1)\right)\,.
\end{equation*}
\end{enumerate}
\end{teora}

The first part of Theorem~\ref{shape_theorem}, accounting for the \emph{quenched
Lyapunov functions}~$\alpha_\lambda$, is taken from \cite{Zer}. We will not
repeat the proof, which relies on the subadditive ergodic theorem, but refer
the reader to the original paper. The second part of the theorem on the
\emph{annealed Lyapunov functions} $\beta_\lambda$ is proven in
Section~\ref{Shape_theorem} with the help of the subadditive limit theorem.

The  Lyapunov functions play an important role in the large deviation
principles. 
For $x\in\mathbb R^d$, we set  
\begin{align*}
I(x)\defi\sup_{\lambda\geq 0}\left(\alpha_\lambda(x)-\lambda\right)
\quad&\text{and}\quad
J(x)\defi\sup_{\lambda\geq 0}\left(\beta_\lambda(x)-\lambda\right)\,.
\intertext{Both the  functions $I$ and $J$ are
continuous and convex increasing on their effective domains}
D_{I}\defi\{x\in\mathbb R^d:\,I(x)<\infty\}
\quad&\text{and}\quad 
D_{J}\defi\{x\in\mathbb R^d:\,J(x)<\infty\}\,,
\end{align*}
of which both equal the closed
unit ball of the $1$-norm in $\mathbb Z^d$ (see p. 272 in \cite{Zer} and
Section~\ref{Large_deviation_principles} of the present article).  In
particular,
$I$ and
$J$ are lower semicontinuous functions with  compact level-sets, which makes
them
\emph{good rate functions} \cite[see e.g.][]{DS}.

\begin{teora}[Large deviation principles]\label{LDP}\noindent
\begin{enumerate}
\item[a)]
There is a set $\Omega'$ of full $\mathbb P$-measure such that for all
$\omega\in\Omega'$ and any drift
$h\in\mathbb R^d$, we have
\begin{equation}\label{limit_qu}
\lim_{n\to\infty}\frac{\log Z^h_{n,\omega}}{n}
=\sup_{x\in\mathbb R^d}\left(h\cdot x-I(x)\right)\,,
\end{equation}
and $S(n)/n$ satisfies a large deviation principle
under $Q^h_{n,\omega}$ with rate~$n$ and good rate function 
\begin{align*}
I_h(x)&\defi I(x)-h\cdot x+\sup_{y\in\mathbb R^d}\left(h\cdot
y-I(y)\right)\,,\quad x\in\mathbb R^d\,,
\end{align*}
as $n$ tends to infinity.
Namely, for any $\omega\in\Omega'$ and  $h\in\mathbb R^d$,
\begin{align*}
\limsup_{n\to\infty}\frac{1}{n}\log Q^h_{n,\omega}[\,S(n)\in n A\,]
&\leq -\inf_{x\in A}I_h(x)\,,\\
\liminf_{n\to\infty}\frac{1}{n}\log Q^h_{n,\omega}[\,S(n)\in n O\,]
&\geq -\inf_{x\in O}I_h(x)
\end{align*}
for all closed
subsets $A\subset\mathbb R^d$ and all open subsets $O\subset\mathbb R^d$.
\newline

\item[b)]
For any drift $h\in\mathbb R^d$, we have
\begin{equation}\label{limit_an}
\lim_{n\to\infty}\frac{\log Z^h_n}{n}
=\sup_{x\in\mathbb R^d}\left(h\cdot x-J(x)\right)\,,
\end{equation}
and $S(n)/n$ satisfies a large deviation principle under
$Q^h_n$  with rate~$n$ and good rate function
\begin{align*}
J_h(x)&\defi J(x)-h\cdot x+\sup_{y\in\mathbb R^d}\left(h\cdot
y-J(y)\right)\,,\quad x\in\mathbb R^d\,,
\end{align*}
as $n$ tends to infinity.
Namely, for any  $h\in\mathbb R^d$,
\begin{align*}
\limsup_{n\to\infty}\frac{1}{n}\log Q^h_n[\,S(n)\in n A\,]
&\leq -\inf_{x\in A}J_h(x)\,,\\
\liminf_{n\to\infty}\frac{1}{n}\log Q^h_n[\,S(n)\in n O\,]
&\geq -\inf_{x\in O}J_h(x)
\end{align*}
for all closed
subsets $A\subset\mathbb R^d$ and all open subsets $O\subset\mathbb R^d$.
\end{enumerate}
\end{teora}

The crucial case of Theorem~\ref{LDP} is the one of  vanishing
drift,  which  for the quenched setting already is proved in \cite{Zer}. 
The extension
to arbitrary  drifts then follows by general principles (essentially 
Varadhan's Lemma).

To describe the  transition between  small and  large drift, we 
quantify the size of $h$ in terms of the dual norms of the Lyapunov functions.
For $\lambda\geq 0$, the  dual norm  of
$\alpha_\lambda$ is defined by 
\begin{align*}
\alpha^*_\lambda(\ell)
&\defi\sup_{x\neq 0}\left(\frac{\ell\cdot x}{\alpha_\lambda(x)}\right)\,,\quad 
\text{$\ell\in\mathbb R^d$}\,,\\
\intertext{while the dual norm  of $\beta_\lambda$ is
defined by}
\beta^*_\lambda(\ell)
&\defi\sup_{x\neq 0}\left(\frac{\ell\cdot x}{\beta_\lambda(x)}\right)\,,\quad 
\text{$\ell\in\mathbb R^d$}\,.
\end{align*}
It is plain to see that  $\alpha^*_\lambda$ and $\beta^*_\lambda$ indeed
are norms  again. Further elementary properties are established in
Section~\ref{dual_norms_phase_transitions}.

As a corollary to Theorem~\ref{shape_theorem}, we have the following ``point to
hyperplane'' interpretation on the dual norms\,: 
For $\ell\neq 0$ and $u\geq  0$, let
$$H_{\ell}(u)\defi\inf\left\{n\geq 0:\, \ell\cdot S(n)\geq u\right\}$$
be  the time of the random walk's first entrance    into
the  half-space $\{x\in\mathbb R^d:\, \ell\cdot x\geq u\}$.

\begin{corola}[Point to hyperplane characterization of dual norms]\label{point_to_hyperplane}
\noindent
\begin{enumerate}
\item[a)] 
There is a set of full $\mathbb P$-measure, on which for all
$\lambda\in[0,\infty)$ and $\ell\in\mathbb R^d\setminus\{0\}$, we have
\begin{equation*}
\lim_{u\to\infty}\frac{1}{u}
\log E\left[\exp\big(-\lambda H_\ell(u)
-\Psi(H_\ell(u),\omega)\big)\right]
=-\frac{1}{\alpha^*_\lambda(\ell)}\,.
\end{equation*}

\item[b)]
For all $\lambda\in[0,\infty)$ and $\ell\in\mathbb R^d\setminus\{0\}$, we have
\begin{equation*}
\lim_{u\to\infty}\frac{1}{u}
\log E\left[\exp\big(-\lambda H_\ell(u)
-\Phi(H_\ell(u))\big)\right] =-\frac{1}{\beta^*_\lambda(\ell)}\,.
\end{equation*}
\end{enumerate}
\end{corola}
Corollary~\ref{point_to_hyperplane}  is the discrete counterpart  to
Sznitman's results for Brownian motion in a Poissonian potential
\cite[Corollary~2.11 and Corollary~3.6. of Chapter~5]{Szn}.

As the following theorem shows, 
the phase transition in the long-time behavior of the random
walk is appropriately characterized by the size of the drift,
measured in terms of  the dual norms $\alpha^*_0$ and $\beta^*_0$.

\begin{teora}[Phase transitions]\label{phase_transitions}
\noindent
\begin{enumerate}
\item[a)]
On the set $\Omega'$ appearing in 
Theorem~\ref{LDP} and for any
$h\in\mathbb R^d$, we have
\begin{equation}\label{supremum_qu}
\lim_{n\to\infty}\frac{1}{n}\log Z^h_{n,\omega}=
\left\{\begin{array}{cl}
0& \text{if $\alpha^*_0(h)\leq 1$}\,,\\
\lambda^{\mathrm{qu}}_h& \text{if $\alpha^*_0(h)> 1$}\,,
\end{array}\right.
\end{equation}
where $\lambda^{\mathrm{qu}}_h>0$  is the unique number with
$\alpha^{*}_{\lambda^{\mathrm{qu}}_h}(h)=1$.
Again on  $\Omega'$, we furthermore have the following limiting behavior\,:
When
$\alpha^*_0(h)<1$, then
\begin{equation*}
\frac{S(n)}{n}\to 0\quad
\text{in $Q^h_{n,\omega}$ probability, as $n\to\infty$}\,.
\end{equation*}
When $\alpha^*_0(h)>1$, then
\begin{equation*}
\mathrm{dist}\left(\frac{ S(n)}{n},\,M_h\right)\to 0\quad
\text{in $Q^h_{n,\omega}$ probability, as $n\to\infty$}\,,
\end{equation*}
where $M_h\defi\{\,x\in\mathbb R^d:\,
h\cdot x-I(x)=\lambda^{\mathrm{qu}}_h\,\}$ is a compact set, which does not
contain the origin.
\newline

\item[b)]
For any $h\in\mathbb R^d$, we have
\begin{align*}
\lim_{n\to\infty}\frac{1}{n}\log Z^h_n=
\left\{\begin{array}{cl}
0& \text{if $\beta^*_0(h)\leq 1$}\,,\\
\lambda^{\mathrm{an\vphantom{qu}}}_h& \text{if $\beta^*_0(h)> 1$}\,,
\end{array}\right.
\end{align*}
where $\lambda^{\mathrm{an\vphantom{qu}}}_h>0$  is the unique 
number with $\beta^{*}_{\lambda^{\mathrm{an\vphantom{qu}}}_h}(h)=1$.
We furthermore have the following limiting behavior:
When $\beta^*_0(h)<1$, 
\begin{equation*}
\frac{ S(n)}{n}\to 0\quad
\text{in $Q^h_n$ probability as $n\to\infty$}\,.
\end{equation*}
When $\beta^*_0(h)>1$, 
\begin{equation*}
\mathrm{dist}\left(\frac{ S(n)}{n},\,N_h\right)\to 0\quad
\text{in $Q^h_n$ probability as $n\to\infty$}\,,
\end{equation*}
where $N_h\defi\{\,x\in\mathbb R^d:\,
h\cdot x-J(x)=\lambda^{\mathrm{an\vphantom{qu}}}_h\,\}$ is a compact set,
which does not contain the origin.
\end{enumerate}
\end{teora}

Remark that for large drifts, since
$M_h$ and $N_h$ are bounded away from the origin, 
Theorem~\ref{phase_transitions}  implies that the random walk $S(n)$  typically
moves away from the origin, with distance  of order $O(n)$ as
$n\to\infty$.  For small drifts, on the other hand, the dislocation rate
$\|S(n)/n\|$ typically falls below any  positive value in the limit
$n\to\infty$.  Theorem~\ref{phase_transitions} thus displays two phase
transitions, in both quenched and annealed settings, between 
\emph{ballistic} behavior of the walk for large $h$ and \emph{sub-ballistic}
behavior for small $h$. 
The unit spheres of
$\alpha^*_0$ and $\beta^*_0$ correspond to the sets of critical drifts.

The normalization for the asymptotics in Theorem~\ref{phase_transitions} is appropriate in the ballistic regime.
In the continuous model, more exact asymptotics for the sub-ballistic phase are established in \cite{Sz95a} and \cite{Sz95b}.
In the  discrete setting, by analogy to the continuous model, we thus believe that convenient normalizations for $\log Z^h_{n,\omega}$, respectively $\log Z^h_n$, are given by $n(\log n)^{-2/d}$, respectively $n^{d/d+2}$.

To conclude this section, we  stress that Theorem~\ref{shape_theorem}-\ref{phase_transitions}  
essentially are  discrete counterparts to Sznitman's results for the Brownian
motion in Poissonian potentials. 
We  however  would like to point out the introduction of $\lambda^{\mathrm{qu}}_h$ and
$\lambda^{\mathrm{an\vphantom{qu}}}_h$ in
Theorem~\ref{phase_transitions},  which we
believe to be new\,: In order to obtain the
ballistic behavior of the random walk, for either the continuous or the
discrete setting, it actually suffices to show that the so-called \emph{Lyapunov
exponents}
$$
\lim_{n\to\infty}\frac{1}{n}\log Z^h_{n,\omega}
\,
\quad\text{and}\quad
\lim_{n\to\infty}\frac{1}{n}\log Z^h_n
$$ are strictly positive (in fact, this is Sznitman's approach).
By means of
$\lambda^{\mathrm{qu}}_h$ and
$\lambda^{\mathrm{an\vphantom{qu}}}_h$, on the other hand, we are able to
express these limits in an implicit way, providing a useful relation to their
counterparts  in the simpler ``point to hyperplane'' setting of
Corollary~\ref{point_to_hyperplane}.  
In fact, in the second, forthcoming paper \cite{Flu}, this relation is used in the context of a renewal formalism to transfer an exponential gap result from  the ``point to hyperplane'' to the ``fixed number of steps'' setting,
implying analyticity of the annealed Lyapunov exponent, and providing  coincidence of the quenched and the annealed exponent for weak potentials in dimensions $d\geq 4$.

\section{Lyapunov functions and shape theorem}
\label{Shape_theorem}

The quenched part of Theorem~\ref{shape_theorem} has been proved by
\cite{Zer}\,: the existence of the norms $\alpha_\lambda$ and the
bounds in
\eqref{bounds_alpha} are  part of Proposition~4, the asymptotic equivalence
in
\eqref{quenched_shape_limit} corresponds to his Theorem~8, and the further
properties of
$\alpha_\lambda$ are established on page 272. Observe that Zerner left out the
condition
$\mathrm{ess\,inf}V_x=0$ instead of introducing the parameter $\lambda$.

In the rest of this section,
we follow Zerner's line to prove the remaining  annealed part of
Theorem~\ref{shape_theorem}.
Recall the two-point function
\begin{equation}\label{b_lambda_section_2}
b_\lambda(x)=
-\log E\left[\exp\big(-\lambda
H(x)-\Phi(H(x))\big);\, H(x)<\infty\right]
\end{equation}
for $\lambda\geq 0$ and  $x\in\mathbb Z^d$. 
The stopping time $H(x)$
denotes the time of the random walk's first visit to the lattice site
$x$,
and the path potential
$\Phi$  is given by
\begin{equation}\label{Phi}
\Phi(n)=\sum_{z\in\mathbb Z^d}\varphi(l_z(n))
\end{equation}
for $n\in\mathbb N$.
Here, $\varphi:[0,\infty)\to[0,\infty)$ is  a non-constant, 
concave increasing function, satisfying 
$\varphi(0)=0$ and
$\lim_{t\to\infty}\varphi(t)/t=0$. 
By dominated convergence and  H\"older inequality, it is plain
that
$b_\lambda(x)$, to any fixed $x\in\mathbb Z^d$,  is  continuous and concave
increasing  in the variable $\lambda\in[0,\infty)$. 
Moreover, we have $H(x)\geq\|x\|_1$ and thus
$\Phi(H(x))\geq\varphi(1)\|x\|_1$ for all  $x\in\mathbb Z^d$.
For any $\lambda\geq 0$, this yields the lower bound
\begin{align}\label{lower_bound_b}
{b_\lambda(x)}
&\geq{\|x\|_1}\left(\lambda+\varphi(1)\right)\,,
\end{align}
while the upper bound
\begin{align}
{b_\lambda(x)}\label{upper_bound_b}
&\leq{\|x\|_1}\left(\log(2 d)+\lambda+\varphi(1)\right)
\end{align}
comes from restricting the
expectation in \eqref{b_lambda_section_2} to a single $\|x\|_1$-step path from
the origin to $x\in\mathbb Z^d$.

To a fixed $\lambda$, we want to establish the triangle inequality for
$b_\lambda$ as a function on $\mathbb Z^d$. 
To this end, let
$$
H(x,z)\defi\inf\left\{m\geq H(x):\,S(m)=z\right\}$$
be the time of the  random walk's first visit to the site
$z\in\mathbb Z^d$ after its first visit to the site $x\in\mathbb Z^d$,
and set
$$\Phi(n,m)\defi\sum_{z\in\mathbb Z^d}\varphi\left(l_z(n)-l_z(m)\right)$$
for $n, m\in\mathbb N_0$ with $n\geq m$.
Again by the concavity of $\varphi$, we  have 
\begin{equation}\label{Phi_inequality}
\Phi(n)\leq \Phi(m)+\Phi(m,n)
\end{equation}
for all $m\geq n$.
The strong Markov
property, applied to the stopping time $H(x)$, then implies
\begin{align}
b_\lambda(x+y)
&\leq
-\log E\left[\,1_{\{H(x)\leq H(x,\,x+y)<\infty\}}\,\exp\big(-\lambda
H(x)-\Phi(H(x))\big)\right.\notag\\
&\hspace{3.8 cm}
\left.
\exp\big(\lambda
H(x,x+y)-\Phi(H(x,x+y))\big)\,\right]\quad\notag\\
&=\,b_\lambda(x)+b_\lambda(y)\,.\label{triangle_inequality}
\end{align}

Given the validity of the triangle inequality \eqref{triangle_inequality},
we can apply the subadditive limit theorem 
\cite[see e.g.][Appendix~II]{Gri}, which guarantees the existence of
 a function $\beta_\lambda:\mathbb Z^d\to[0,\infty)$
such that
\begin{equation}\label{subadditive_limit}
\lim_{n\to\infty} \frac{1}{n}\,b_\lambda(n x)
=\inf_{n\in\mathbb N} \frac{1}{n}\,b_\lambda(n x)
=\beta_\lambda(x)
\end{equation}
for every $x\in\mathbb Z^d$.
It is easy to conclude that
$\beta_\lambda$ inherits from
$b_\lambda$ the same bounds as in
\eqref{lower_bound_b} and \eqref{upper_bound_b}, that is
\begin{equation}\label{bounds_beta_2}
\lambda+\varphi(1)\leq\frac{\beta_\lambda(x)}{\|x\|_1}
\leq\log(2 d)+\lambda+\varphi(1)
\end{equation}
for all $x\in\mathbb Z^d\setminus\{0\}$, and that
\begin{equation}\label{properties_beta_2}
\begin{split}
\beta_\lambda(n x)&=n\beta_\lambda(x)\,,\\ 
\beta_\lambda(x+y)&\leq\beta_\lambda(x)+\beta_\lambda(y)
\end{split}
\end{equation}
are satisfied for any $n\in\mathbb N$ and $x,y\in\mathbb Z$.
Moreover, to fixed $x\in\mathbb Z^d$, $\beta_\lambda(x)$  is continuous and
concave increasing in $\lambda\in[0,\infty)$\,:
As a limit of
concave functions,
$\beta_\lambda(x)$ is concave again and thus lower semicontinuous (possibly being discontinuous in $\lambda=0$). The
upper semicontinuity, by the representation of $\beta_\lambda(x)$ as an infimum
in \eqref{subadditive_limit}, is derived from the continuity of $b_\lambda(n x)$
in
$\lambda$ for $n\in\mathbb N$.

By setting $\beta_\lambda(q x)=q\beta_\lambda(x)$ for $q\in\mathbb Q$, we
extend
$\beta_\lambda$ well-defined at first to a function on $\mathbb Q^d$ and then
by continuity to a function on $\mathbb R^d$. Thereby, $\beta_\lambda$
maintains its properties as a function of $\lambda$ and still satisfies
\eqref{bounds_beta_2} and \eqref{properties_beta_2}. 
In particular,
$\beta_\lambda$ is a norm on $\mathbb R^d$.
Moreover, since
$$\left|\beta_{\lambda_k}(x_k)-\beta_{\lambda}(x)\right|
\leq
\left|\beta_{\max_{k\in\mathbb N}\lambda_k}(x_k-x)\right|
+\left|\beta_{\lambda_k}(x)-\beta_{\lambda}(x)\right|$$
for all sequences $(\lambda_k)_{k\in\mathbb N}$  and $(x_k)_{k\in\mathbb N}$
with $\lambda_k\to\lambda$ and  $x_k\to x$,  we obtain the joint continuity
of $\beta_\lambda(x)$  in
$(\lambda,x)\in[0,\infty)\times\mathbb R^d$ from the
continuity  in the single arguments.

It remains to prove the limiting behavior of $b_\lambda/\beta_\lambda$ in 
\eqref{annealed_shape_limit}. It suffices to show 
\begin{equation}\label{sufficient_formulation}
\lim_{k\to\infty}
\left|\frac{b_\lambda(x_k)-\beta_\lambda(x_k)}{\|x_k\|_1}\right|=0
\end{equation}
for any  sequence $(x_k)_{k\in\mathbb N}$ on $\mathbb Z^d$ with
$\|x_k\|_1\to\infty$. We yet can restrict to the case where
$x_k/\|x_k\|_1\to e$ for some point $e\in S^{d-1}$; if
\eqref{sufficient_formulation} was not true for an arbitrary sequence, it
would not be true
 for a subsequence with this convergence property either.

To this end, for any $\varepsilon>0$, choose $\tilde e\in\mathbb Q^d$ and
$m\in\mathbb N$ such that $m \tilde e\in\mathbb Z^d$  and
$\|e-\tilde e\|_1<\varepsilon$ as well as
$|\beta_\lambda(e)-\beta_\lambda(\tilde e)|<\varepsilon$. We approximate
$(x_k)_{k\in\mathbb N}$ by the sequence
$(n_k\,x)_{k\in\mathbb N}$ on $\mathbb Z^d$, where
$$x=m \tilde e\in\mathbb Z^d\quad\text{and}\quad
n_k=\left\lfloor\frac{\|x_k\|_1}{m}\right\rfloor\,.$$ Thereby
$\lfloor\cdot\rfloor$ denotes the largest integer less than or equal to a
real number. 
Notice first that  $\lim_{k\to\infty}\|x_k\|_1/n_k=m$. 
We thus have
\begin{align}
\|x_k-n_k\,x\|_1
&\leq
\left\|x_k- \frac{n_k\,m}{\|x_k\|_1}\,x_k\right\|_1
+\left\|
\frac{n_k\,m}{\|x_k\|_1}\,x_k-n_k\,x\right\|_1 \notag\\
&=
\left(1-\frac{n_k\,m}{\|x_k\|_1}\right)\|x_k\|_1
+ n_k\,m
\left\|\frac{x_k}{\|x_k\|_1}-\tilde e\right\|_1\notag\\
&<\,\varepsilon \|x_k\|_1\label{smaller}
\end{align}
for $k$ large enough.
By 
the (inverted) triangle inequality \eqref{triangle_inequality} for $b_\lambda$,
we get
\begin{multline*}
\left|\frac{b_\lambda(x_k)-\beta_\lambda(x_k)}{\|x_k\|_1}\right|
\leq
\,\frac{b_\lambda(x_k-n_k\,x)}{\|x_k\|_1}\\
\,+\left|\frac{b_\lambda(n_k\,x)}{\|x_k\|_1}-\beta_\lambda(\tilde e)\right|
+\left|\beta_\lambda(\tilde e)
-\beta_\lambda\left(\frac{x_k}{\|x_k\|_1}\right)\right|\,.
\end{multline*}
The first summand on the right-hand side is bounded from above by $\varepsilon
c_\lambda$ with  $c_\lambda=\log(2 d)+\lambda+\varphi(1)$ due to
\eqref{upper_bound_b} and
\eqref{smaller}.  The second summand  tends to zero for $k$ going to infinity
since
$\|x_k\|_1/n_k\to m$ and
$b_\lambda(n_k\,x)/n_k\to \beta_\lambda( x)$.
The last summand finally is smaller than $\varepsilon$ for $k$ large enough
by the  assumptions  $\lim_{k\to\infty}x_k/\|x_k\|_1=e$ and
$|\beta_\lambda(\tilde e)-\beta_\lambda(e)|<\varepsilon$. Hence, letting
$\varepsilon$ tend to zero implies \eqref{sufficient_formulation} and completes
the proof of the shape theorem in the annealed setting.

\section{Large deviation principles}\label{Large_deviation_principles}

The aim of this section is to prove Theorem~\ref{LDP}.
The limit results \eqref{limit_qu} and
\eqref{limit_an} for arbitrary drifts as well as the LDP's for
$Q^{h\neq 0}_n$ and $Q^{h\neq 0}_{n,\omega}$ thereby follow from the LDP's for
$Q^{h= 0}_n$ and $Q^{h= 0}_{n,\omega}$ as an application of Varadhan's
lemma \citep[see e.g.][Theorem 2.1.10 and Exercise 2.1.24]{DS}. To this
purpose, we only need to establish the  ``exponential tightness estimates''
\begin{align*}
\lim_{L\to\infty}\limsup_{n\to\infty}\frac{1}{n}\log
E\left[\exp\big(h\cdot S(n)-\Psi(n,\omega)\big);\,h\cdot S(n)\geq n
L\right]=-\infty\,,
\end{align*}
for the quenched setting and
\begin{align*}
\lim_{L\to\infty}\limsup_{n\to\infty}\frac{1}{n}\log
E\left[\exp\big(h\cdot S(n)- \Phi(n)\big);\,
h\cdot S(n)\geq n L\right]=-\infty
\end{align*}
for the annealed setting .
But, since both the expectations in the above limits  are bounded
by
\begin{align*}
E\,\big[\exp(h\cdot S(n));\,h\cdot S(n)\geq n L\big]
&\leq \exp(-n L)  E\,\big[\exp(2 h\cdot S(n))\big]\\
&=  \exp(-n L) {E \big[\exp(2 h\cdot S(1))\big]}^{n}\,,
\end{align*}
the exponential estimates follow immediately.

For vanishing drift, the limit in \eqref{limit_qu} and
the large deviation property in the quenched setting have already been
proved
\citep[Proposition~17 and Theorem~19]{Zer}.  

We follow Zerner's line  to prove
the remaining annealed part of Theorem~\ref{LDP} for the case $h=0$.
That is, we investigate the large deviations of the symmetric
random walk under the annealed path measures $Q_n$ with density
\begin{align*}
\frac{d Q_n}{d P}
\defi\frac{\exp(-\Phi(n))}{Z_n}
\end{align*}
{when $n\in\mathbb N$ tends to infinity, 
where the  normalization constant $Z_n$  is given by}
\begin{align*}
Z_n \defi E\,\big[\exp(-\Phi(n))\big]
=E\left[\exp\big(-\sum\nolimits_{x\in\mathbb Z^d}\varphi(l_x(n)\big)\right]\,.
\end{align*}
Thereby, we
have
$Q_n=Q^{h=0}_n$ and  $Z_n=Z^{h=0}_n$  according to the notations from
Section~\ref{Introduction}.

We first take care of the normalization constant $Z_n$. 
Claim \eqref{limit_an} in Theorem~\ref{LDP} clearly reduces to
\begin{equation}\label{vanishing}
\lim_{n\to\infty}\frac{-\log
Z_n}{n}=0\,.
\end{equation}
In fact, it turns out that the above limit equals
$\lim_{n\to\infty}\varphi(n)/n$, which is assumed to be zero.
To see this, observe that 
\begin{align*}
\lim_{n\to\infty}\frac{\varphi(n)}{n}&=
\inf_{n\in\mathbb N}\frac{\varphi(n)}{n}\,,
\intertext{once again by the concavity  of $\varphi$.
By the definition of
$\Phi$,  we therefore have}
\lim_{n\to\infty}\frac{-\log Z_n}{n}
&\geq\lim_{n\to\infty}\frac{\varphi(n)}{n}\,.
\end{align*}

It remains to prove the upper estimate.
For any integer $R$ and all $n\in\mathbb N$, we obviously have
\begin{align*}
Z_n\geq E\,\big[\exp(-\Phi(n));\,
\text{$\|S(m)\|_1\leq R$ for $m\leq n$}\big]\,.
\end{align*}
In order to find a lower bound for the right-hand side  of this inequality,
observe that 
$$
\sum_{x\in\mathbb Z^d}\varphi(l_{x}(n))\,
1_{\left\{\text{$\|S(m)\|_1\leq R$ for $m\leq n$}\right\}}
\leq (2 R+1)^d \varphi(n)$$ 
is valid for all $n\in\mathbb N$, and that 
a $2R$-step path  with start and end at the origin remains
within the cube
$\{x\in\mathbb Z^d:\|x\|_1\leq R\}$.  
By  the Markov property, we thus obtain
\begin{align*}
&-\limsup_{n\to\infty}\frac{1}{n}\,
\log E\,\big[\exp(-\Phi(n));\,
\text{$\|S(m)\|_1\leq R$ for $m\leq n$}\big]\\
&\qquad\leq
\lim_{n\to\infty}\frac{\varphi(n)}{n}- 
\limsup_{n\to\infty}\frac{1}{n}\,\log P\big[\,\text{$\|S(m)\|_1\leq R$ for
$m\leq n$}\big]\\
&\qquad\leq
\lim_{n\to\infty}\frac{\varphi(n)}{n}-
\frac{1}{2R}\, \log P[\,S(2R)=0\,]\,,
\end{align*}
of which the last summand vanishes when $R$ tends to infinity by the local
central limit theorem \cite[see e.g.][]{Woe}. This completes the
proof of~\eqref{vanishing}.

We step forward to the  large deviation
principle.
For the case $h=0$, the rate
function  will be 
$$J(x)=\sup_{\lambda\geq 0}\left(\beta_{\lambda}(x)-\lambda\right)
\,,\quad x\in\mathbb R^d\,.$$ 
As a supremum of continuous functions, $J$ is lower
semicontinuous.  Furthermore, $J$ inherits the
convexity  from the norms
$\beta_\lambda$ and hence is upper semicontinuous on
its effective
domain 
$$D_J=\{x\in\mathbb R^d:\, J(x)<\infty\}\,.$$  
Moreover, the bounds for $\beta_\lambda$ in \eqref{bounds_beta_2} yield
that
$D_J$ equals the closed unit ball of the $1$-norm.

The rest of this section is devoted to the proof of the large
deviations estimates\,:
For any closed
subset $A\subset\mathbb R^d$ and  open subset $O\subset\mathbb R^d$,
\begin{align}\label{upper_bound}
\limsup_{n\to\infty}\frac{1}{n}\log Q_n[\,S(n)\in n A\,]&\leq -\inf_{x\in
A}J(x)\,,\\
\label{lower_bound}
\liminf_{n\to\infty}\frac{1}{n}\log Q_n[\,S(n)\in n O\,]&\geq -\inf_{x\in
O}J(x)\,.
\end{align}

We start with the upper estimate.
Since $J(x)=\infty$ if $\|x\|>1$ and $\|S(n)\|_1\leq n$ for $n\in\mathbb
N$, we can restrict to the case where $A\subset D_J$ is compact.
For $n\in\mathbb N$, we set $H(n
A)=\inf\{H(x):\, x\in n A\}$.
 Since $\{S(n)\in n A\}\subset\{H(n A)\leq n\}$, we have
\begin{equation*}
\log E\,\big[\exp(-\Phi(n));\, S(n)\in n A\big]
\leq -(b_\lambda(n A)-\lambda n)
\end{equation*}
for all $\lambda\geq 0$, where $b_\lambda(n A)$ is defined as in
\eqref{b_lambda_section_2}, but with $H(y)$ replaced by $H(n A)$.  
From the representation in
\eqref{subadditive_limit}  of
$\beta_\lambda(x)$ as an infimum, and since $A$ is bounded, we obtain
\begin{align*}
\frac{b_\lambda(n A)}{n}
&\geq \frac{-\log |n A\cap \mathbb Z^d|}{n} +\sup_{x\in  A\cap
\frac{1}{n}\mathbb Z^d}\, \frac{b_\lambda(n x)}{n}\\ 
&\geq \frac{-\log (n K)}{n} +\sup_{x\in A}
\beta_\lambda(a)
\end{align*} 
for all $n\in\mathbb N$ and some constant $K=K(A,d)$.
By \eqref{vanishing}, we then have 
\begin{align}\label{upper_1}
\limsup_{n\to\infty}\frac{1}{n}\,
\log Q_n[\,S(n)\in n A\,]
&=-\sup_{\lambda\geq 0}\,\inf_{x\in A}\beta_\lambda(x)-\lambda\,.
\end{align}

However, in order to complete the proof of \eqref{upper_bound}, we  need to
exchange infimum and supremum in \eqref{upper_1}.
For any $\varepsilon>0$, thanks to the compactness of $A$, there are
$m\in\mathbb N$ and $\lambda_1,\dots,\lambda_m>0$ such that the compact sets
$$
A_i=\big\{y\in A:\,\beta_{\lambda_i}(y)-\lambda_i\geq \inf_{x\in A}
J(x)-\varepsilon\big\}\,,\quad i=1,\dots,m\,,$$
cover $A$.
From \eqref{upper_1} applied to the sets $A_i$, we therefore obtain 
\begin{align*}
\limsup_{n\to\infty}\frac{1}{n}\,
\log Q_n[\,S(n)\in n A\,]
&\leq  \max_{i=1,\dots,m}\, \limsup_{n\to\infty}\frac{1}{n}\,
\log Q_n[\,S(n)\in n A_i\,]\\
&\leq -\min_{i=1,\dots,m}\,
\inf_{x\in A_i}\beta_{\lambda_i}(x)-\lambda_i\\
&\leq-\inf_{y\in A}J(x)+\varepsilon\,.
\end{align*}
Since $\varepsilon>0$ was arbitrary, this  proves \eqref{upper_bound}.

We step forward to the proof of \eqref{lower_bound}. That is, for an open set
$O\subset\mathbb R^d$, we need to show
\begin{equation}\label{to_show}
\liminf_{n\to\infty}\frac{1}{n}\log Q_n[\,S_n\in n O\,]\geq
-J(x)
\end{equation}
for all $x\in O\cap D_J\setminus\{0\}$,
where the origin  can be  excluded since   $J$ is
continuous on $D_J$. To do this, we
first determine some number $\lambda_0$, at which $\beta_\lambda(x)-\lambda$
attains its maximum as a function of $\lambda$\,:
Since $\beta_\lambda(x)$ is concave in $\lambda$,
the right derivative 
$$ 
\dot\beta^+_\lambda(x)=\lim_{\varepsilon\downarrow 0}\,
\frac{\beta_{\lambda+\varepsilon}(x)-\beta_{\lambda}(x)}{\varepsilon}
$$  
of $\beta_\lambda(x)$
is a well-defined and nondecreasing (but not necessarily continuous) function
of
$\lambda\in[0,\infty)$. 
Now, if 
$\dot\beta^+_\lambda(x)<1$ for all
$\lambda> 0$, the  maximum is located at $\lambda_0=0$.
Otherwise we choose
$$
\lambda_0=\inf\,\{\lambda>
0:\,\dot\beta^+_{\lambda}(x)< 1\}\,,
$$
which is finite since $x\in D_J\cap O$,
and which is a transition point from nondecreasing
to decreasing behavior for the map
$\lambda\mapsto\beta_\lambda(x)-\lambda$. Hence, in both cases, we have
\begin{equation}\label{lambda_0} 
I(x)=\beta_{\lambda_0}(x)-\lambda_0\,.
\end{equation}

We now take care of the fact that  
$\dot\beta^+_\lambda(x)$ might be discontinuous in $\lambda_0$.
In that case,  there is  a non
trivial  ``interval'' in the  half-line
$\{s x:\,s\in\mathbb [0,\infty)\}$, on which \eqref{lambda_0} remains true
with the fixed constant $\lambda_0$. We  express $x$ by a linear combination
of the end points of this interval. 
Let $\dot\beta^-_\lambda(x)$ denote the
left derivative of
$\beta_\lambda(x)$ with respect to $\lambda>0$.
Set 
$$
y^+=\frac{x}{\dot\beta^+_{\lambda_0}(x)}\quad\text{and}\quad
y^-=
\left\{\begin{array}{cl}
0& \text{if $\lambda_0=0$}\,,\\
{\frac{x\vphantom{\dot\beta^-_{\lambda_0}(x)}}{\dot\beta^-_{\lambda_0}(x)}}&
\text{if $\lambda_0>0$}\,.
\end{array}\right.
$$
We then have $$(1-t)y^-+t y^+=x$$
with $t=\frac{|y-y^-|}{|y^+-y^-|}<1$ if
$\dot\beta^+_{\lambda_0}(x)<\dot\beta^-_{\lambda_0}(x)$, and with 
$t=1$ in the continuous case. 
For a reason that will become clear later, we
furthermore approximate $y^-$ by slightly ``smaller'' sites
$y^-_\rho=\rho y^-$. Since $O$ is open, we can choose
$\rho<1$ large enough to fulfill  
\begin{equation}\label{contained_1}
(1-t)y^-_\rho+t y^+\in O.
\end{equation}
Let finally $(y^-_{\rho,n})$ and $(y^+_n)$ be two sequences in $\mathbb Z^d$
such that
\begin{equation}\label{contained_2}
\lim_{n\to\infty}\frac{y^-_{\rho,n}}{n}=(1-t)y^-_{\rho}\,\quad\text{and}\quad
\lim_{n\to\infty}\frac{y^+_n}{n}=t y^+,
\end{equation}
and set $x_n=y^-_{\rho,n}+ y^+_n$ for $n\in\mathbb N$. 
Thereby, if $\lambda_0=0$
or $t=1$,  we     may simply set
$y^-_{\rho,n}=0$  for
$n\in\mathbb N$.

We want to renew the Markov chain at the sites $y^-_{\rho,n}$ and $x_n$.
To this end, let $R$ be an arbitrary integer. 
Since $O$ is open, we obtain from \eqref{contained_1} and
\eqref{contained_2} that there exists some $n_0\in\mathbb N$ such that
$x_n+y\in n O$ is valid for all $\|y\|_1\leq R$ and all $n\geq n_0$.
As a consequence,
\begin{multline*}
\big\{H(y^-_{\rho,n})\leq(1-t)n\big\}
\cap\big\{H(y^-_{\rho,n},x_n)\leq n\big\}\\
\cap\big\{\|S(m)-x_n\|_1\leq R\text{ for
$H(y^-_{\rho,n},x_n)<m\leq H(y^-_{\rho,n},x_n)+n$}\big\}
\end{multline*}
is contained in $\{S(n)\in n O\}$ for $n$ large enough.
By the monotonicity of $\Phi$ and a double application of
\eqref{Phi_inequality}, we furthermore have
\begin{multline*}
\Phi(n)\leq
\Phi\big(H(y^-_{\rho,n})\big)+
\Phi\big(H(y^-_{\rho,n}),H(y^-_{\rho,n},x_n)\big)\\+
\Phi\big(H(y^-_{\rho,n},x_n),H(y^-_{\rho,n},x_n)+n\big)\,.
\end{multline*}
From \eqref{vanishing} and the strong Markov property, it 
thus follows that the left-hand side of \eqref{lower_bound} is not smaller
than
\begin{align}\label{split_1}
&\liminf_{n\to\infty}\frac{1}{n}\,
\log E\left[\exp\left(-\Phi(H(y^-_{\rho,n}))\right);\,
H(y^-_{\rho,n})\leq (1-t)n\right]\\ 
&\qquad\quad\label{split_2}
+\liminf_{n\to\infty}\frac{1}{n}\,
\log E\left[\exp\left(-\Phi(H(y^+_{n}))\right);\, H(y^+_{n})\leq
t n\right]\\ 
&\qquad\quad\label{split_3}
+\liminf_{n\to\infty}\frac{1}{n}\,
\log E\Big[\exp\big(-\Phi(n)\big);\,
\text{$\|S(m)\|_1\leq R$ for $m\leq n$}\Big]\,,\notag
\end{align}
of which the last summand vanishes when $R$ tends to infinity, as we have seen
in the proof of \eqref{vanishing}.

In order to bound the first and the second summand, we need the following
result\,: For $\lambda\geq 0$ and $y\in\mathbb Z^d$, let the distribution
$P^y_\lambda$ be given by means of the density
$$\frac{d P^y_\lambda}{d P}
=\frac{1}{Z^y_\lambda}\,\exp\big(-\lambda H(y)-\Phi(H(y))\big)\,
1_{\{H(y)<\infty\}}\,,$$
where $Z^y_\lambda=b_\lambda(y)$ is the corresponding normalization
constant.

\begin{lemma}\label{CLT}
Let  $(y_n)$ a be a sequence in $\mathbb Z^d$ with
$\lim_{n\to\infty}{y_n}/{n}=y\in\mathbb R^d\setminus\{0\}$. 
For every $\lambda>0$ and
$\gamma<\dot\beta^+_{\lambda}(y)\leq
\dot\beta^-_{\lambda}(y)<\delta$, we  have
\begin{equation}\label{large_numbers}
\lim_{n\to\infty} P^{y_n}_\lambda\big[H(y_n)/n\in[\gamma,\delta]\big]=1\,.
\end{equation}
\end{lemma}

{\bf Proof.\,}
For any $0<\varepsilon<\lambda$, we
have
\begin{align*}
&E\left[\exp\big(-\lambda
H(y_n)-\Phi(H(y_n))\big);\,H(y_n)\notin {{[\gamma n,\delta n]}}
\right]\\ 
&\hspace{.7 cm}\leq	
\exp(\varepsilon\gamma n)\,
E\left[\exp\big(-(\lambda+\varepsilon) H(y_n)-\Phi(H(y_n))\big);\,  
H(y_n)<\infty\right]\\
&\hspace{1.2 cm}+
\exp(-\varepsilon \delta n)\,
E\left[\exp\big(-(\lambda-\varepsilon)
H(y_n)-\Phi(H(y_n))\big);\,   H(y_n)<\infty\right]\,.
\end{align*}
From the shape theorem therefore follows
\begin{multline*}
\limsup_{n\to\infty}\frac{1}{n}
\log P^{y_n}_\lambda\big[{H(y_n)}/n\notin
{[\gamma,\delta]}\big]\\
\leq 
-\varepsilon\min\bigg\{
\frac{\beta_{\lambda+\varepsilon}(y)-\beta_{\lambda}(y)}{\varepsilon}-\gamma
\,,\,
\delta-\frac{\beta_{\varepsilon}(y)-\beta_{\lambda-\varepsilon}(y)}{\varepsilon}
\bigg\}\,.
\end{multline*}
By the
assumptions on $\gamma$ and $\delta$, the right-hand side of this last
expression is strictly negative for
$\varepsilon>0$ small enough, which then implies \eqref{large_numbers}.
\hfill{$\qed$}

We are now able to complete the proof of \eqref{lower_bound}.
Suppose  $\lambda> \lambda_0$ and
$\gamma<\dot\beta^+_\lambda(y^+)$.
The shape theorem  provides that \eqref{split_2}
is not smaller than
$$
\gamma\lambda-\beta_\lambda(y^+)+
\lim_{n\to\infty}\frac{1}{n}\log
P^{y^+_n}_\lambda\big[H(y_n)/n\in[\gamma,t]\big]\,,
$$
for which Lemma~\ref{CLT} applies  because of
$\dot\beta^-_\lambda(y^+)<\dot\beta^+_{\lambda_0}(y^+)=t$, the strict
inequality coming from the choice of $\lambda_0$.  Since
$\dot\beta^+_\lambda(y^+)$ is upper semicontinuous in
$\lambda_0$, we  thus obtain that
\eqref{split_2} is not smaller than
\begin{align*}
\sup_{\lambda>\lambda_0\vphantom{\dot\beta^+_\lambda(y^+)}}\,
\sup_{\gamma<\dot\beta^+_\lambda(y^+)}\gamma 
\lambda -\beta_\lambda(t y^+)
&= t\lambda_0-\beta_{\lambda_0}(t y^+)\,.
\end{align*}
If $\lambda_0=0$ or $t=1$ is the case, by setting $y^-_{\rho,n}=0$ for 
$n\in\mathbb N$, this already proves~\eqref{to_show}. 
Suppose now $\lambda_0>0$ and $t<1$, which implies $y^-_\rho\neq 0$. 
Since $\dot\beta^-_\lambda(y^-)$ is lower semicontinuous in $\lambda_0$, we
have
$\rho(1-t)=\dot\beta^-_{\lambda_0}(y^-_\rho)\leq\dot\beta^+_{\lambda}(y^-_\rho)$
and $\dot\beta^-_{\lambda}(y^-_\rho)<(1-t)$
whenever
$\lambda<\lambda_0$ is large enough.
The shape theorem and  Lemma~\ref{CLT} with
$\gamma=\rho^2(1-t)$ and $\delta=(1-t)$ then imply that
\eqref{split_1} is not smaller than  
$$
\sup_{\lambda>\lambda_0}\rho^2(1-t)\lambda_0-\beta_\lambda((1-t)y^-_\rho)
= \rho^2(1-t)\lambda_0-\beta_{\lambda_0}(\rho(1-t)y^-)\,.
$$
Since $\rho<1$ was arbitrary, we obtain \eqref{to_show}. This completes the
proof of \eqref{lower_bound}.

\section{Dual norms and phase transitions}
\label{dual_norms_phase_transitions}

The aim of this section is to prove  
Corollary~\ref{point_to_hyperplane} and Theorem~\ref{phase_transitions}.
Recall the definition of the dual norms
$$\alpha^*_\lambda(\ell)
=\sup_{x\neq 0}\left(\frac{\ell\cdot x}{\alpha_\lambda(x)}\right)
\quad\text{and}\quad
\beta^*_\lambda(\ell)
=\sup_{x\neq 0}\left(\frac{\ell\cdot x}{\beta_\lambda(x)}\right)$$
for $\lambda\geq 0$ and $\ell\in\mathbb R^d$.
We first prove some elementary properties of $\alpha^*_\lambda$ and
$\beta^*_\lambda$, similar to 
the ones  of the Lyapunov functions $\alpha_\lambda$ and
$\beta_\lambda$ in
Theorem~\ref{shape_theorem}.

\begin{lemma}\label{properties_dual_norm}
\begin{enumerate}
\item[a)]
$1/\alpha^*_\lambda(\ell)$ is 
continuous in
$(\lambda,\ell)\in[0,\infty)\times \mathbb R^d$ and concave increasing in
$\lambda\in[0,\infty)$,  satisfying
\begin{equation*}
\frac{\|\ell\|_1}{\lambda+\log(2 d)+\mathbb EV_x}
\leq{\alpha^*_\lambda(\ell)}
\leq
\frac{\|\ell\|_1}{\lambda-\log\mathbb E\,\exp(-V_x)}\,.
\end{equation*}

\item[b)]
$1/\beta^*_\lambda(\ell)$ is 
continuous in
$(\lambda,\ell)\in[0,\infty)\times \mathbb R^d$ and concave increasing in
$\lambda\in[0,\infty)$,  satisfying
\begin{equation*}
\frac{\|\ell\|_1}{\lambda+\log(2 d)+\varphi(1)}
\leq{\beta^*_\lambda(\ell)}
\leq
\frac{\|\ell\|_1}{\lambda+\varphi(1)}\,.
\end{equation*}
\end{enumerate}
\end{lemma}

{\bf Proof.\,}
Since the proof works the same way for either the quenched or the annealed
case, we  can restrict to  the quenched setting.
By the 
definition of $\alpha^{*}_\lambda(\ell)$,
to fixed $\ell\in\mathbb R^d$, the
concavity  of
$\alpha^{*}_\lambda(\ell)$ in $\lambda$ is derived from the concavity of
$\alpha_\lambda(x)$ in $ \lambda$ to every fixed $x\in\mathbb R^d$.  The
concavity then  implies lower semicontinuity in $\lambda$, while the upper
semicontinuity, again by the definition of
$\alpha^{*}_\lambda(\ell)$, is derived from the continuity of
$\alpha_\lambda(x)$ in $\lambda$. This proves continuity in the $\lambda$
variable; continuity in the $x$ variable is obvious. 
Let now $(\lambda_k)_{k\in\mathbb N}$ and $(\ell_k)_{k\in\mathbb N}$ be two
sequences with $\lambda_k\to\lambda$ and $\ell_k\to\ell$. We then have
$$\left|\alpha^*_{\lambda_k}(\ell_k)-\alpha^*_{\lambda}(\ell)\right|
\leq
\left|\alpha^*_{\max_{k\in\mathbb N}\lambda_k}(\ell_k-\ell)\right|
+\left|\alpha^*_{\lambda_k}(\ell)-\alpha^*_{\lambda}(\ell)\right|$$
and thus
$\lim_{k\to\infty}\alpha^*_{\lambda_k}(\ell_k)=\alpha^*_{\lambda}(\ell)$. 
This
proves the joint continuity. 
The
bounds for
$\alpha^*_\lambda(\ell)$ finally follow from the bounds for
$\alpha_\lambda(x)$ in
\eqref{bounds_alpha} by standard calculations.
\hfill{$\qed$}

The ``point to hyperplane'' interpretation on the dual norms
in Corollary~\ref{point_to_hyperplane} is derived from  the shape
theorem (Theorem~\ref{shape_theorem}).  The proof is a
modification of   Sznitman's proof for the continuous setting \citep{Szn}.
Since it works in a similar way for either the quenched or annealed case, we
restrict to the more complex  quenched model.  Here,
we have to find  a set of full $\mathbb P$-measure, on which
\begin{equation}\label{hyperplane_quenched}
\lim_{u\to\infty}\frac{1}{u}
\log E\left[\exp\big(-\lambda H_\ell(u)
-\Psi(H_\ell(u),\omega)\big)\right]
=-\frac{1}{\alpha^*_\lambda(\ell)}
\end{equation}
for all
$\lambda\in[0,\infty)$ and $\ell\in\mathbb R^d\setminus\{0\}$, where 
 $H_{\ell}(u)$ is the  time of first entrance into
the half-space $\{x\in\mathbb R^d:\, \ell\cdot x\geq u\}$.
To this end, let 
$(u_n)_{n\in\mathbb N}$ be an arbitrary sequence of  numbers with
$\lim_{n\to\infty} u_n=\infty$. For any fixed $\lambda\in[0,\infty)$, by  the
scalar linearity of  $\alpha_\lambda$ on $\mathbb R^d$, we have
$$\alpha^*_\lambda(\ell)=
\sup_{x\in\mathbb R^d:\,\ell\cdot x=1}\,
\frac{1}{\alpha_\lambda(x)}\,.$$
Consequently, since $\alpha_\lambda$ is continuous and
$\lim_{\|x\|_1\to\infty}\alpha_\lambda(x)=\infty$, there exists
$x^*\in\mathbb R^d$ such that $\ell\cdot x^*=1$ and
$$\alpha^*_\lambda(\ell)=\frac{1}{\alpha_\lambda(x^*)}\,.$$ 

In order to find a lower bound  for the left-hand side of
\eqref{hyperplane_quenched},
choose a sequence $(x_n)_{n\in\mathbb N}$   in
$\mathbb Z^d$ such that $x_n\cdot \ell \geq u_n$ and $\lim_{n\to\infty}
x_n/u_n=x^*$. We then have $H_\ell(u_n)\leq H(x_n)$  and thus
$$
\log E\left[\exp\big(-\lambda H_\ell(u)
-\Phi(H_\ell(u),\omega)\big)\right]
\geq -a_\lambda(x_n,\omega)$$
for all $\omega\in\Omega$ and $n\in\mathbb N$.
On  the set $\Omega_\lambda$ of full
$\mathbb P$-measure appearing in Theorem~\ref{shape_theorem}, we consequently
have  
$$\lim_{n\to\infty}\frac{a_\lambda(x_n,\omega)}{u_n}=\alpha_\lambda(x^*)\,.$$
Since the set $\Omega_\lambda$ does not depend on the sequence
$(x_n)_{n\in\mathbb N}$,  this proves the lower bound part of
\eqref{hyperplane_quenched} on $\Omega_\lambda$ for a fixed $\lambda$ and all
$\ell\neq 0$.

For the upper estimate of the left-hand side  of \eqref{hyperplane_quenched},
choose a number $R$ large enough such that
\begin{equation}\label{condition_R}
\inf\left\{\alpha_\lambda(x):\,\|x\|_1\geq R\right\}
\geq \alpha_\lambda(x^*)\,,
\end{equation}
which is possible since
$\alpha_\lambda$ is a norm. For $\omega\in\Omega$ and any $K\subset\mathbb Z^d$,
we set
\begin{align*}
a_\lambda(K,\omega)=-\log E\left[\exp\big(-\lambda
H(K)-\Psi(H(K),\omega)\big);\, H(K)<\infty\right]
\end{align*}
with $H(K)=\inf\{H(y):\,y\in K\}$. 
 For $n\in\mathbb N$, we furthermore set
$$D_n=\{x\in\mathbb Z^d:\, \|x\|_1\geq R u_n\}
\cap\{x\in\mathbb Z^d:\, \ell\cdot x \geq u_n\}\,,$$ 
whose interior boundary is
$$L_n=\{x\in D_n:\,\text{$\|x-y\|_1=1$ for some $y\neq D_n$}\}\,.$$
Since we have $H(D_n)\leq H_\ell(u_n)$, the logarithmic expectation in 
\eqref{hyperplane_quenched} is bounded from above by
\begin{equation}\label{f}
-a_\lambda(D_n,\omega)=-a_\lambda(L_n,\omega)
\leq \log|L_n|-\min_{x\in L_n}a_\lambda(x,\omega)\,.
\end{equation}
Now, since $|L_n|\leq (2 R u_n+1)^d$, it only remains to take care of
the minimum in \eqref{f}, which we assume to be attained at a site $x_n\in
L_n$.  Again by the shape theorem, we have
$$
\liminf_{n\to\infty}\frac{a_\lambda(x_n,\omega)}{u_n}
=\liminf_{n\to\infty}\frac{\alpha(x_n)}{u_n}
\geq\liminf_{n\to\infty}\inf_{x\in D_n}\frac{\alpha_\lambda(x)}{u_n}
\geq \alpha_\lambda(x^*)$$
on the same $\Omega_\lambda$ of full $\mathbb P$-measure as before,
where the last estimate follows from  \eqref{condition_R}. 
This completes the proof of
\eqref{hyperplane_quenched} on $\Omega_\lambda$ for a fixed $\lambda$ and all
$\ell\neq 0$.

It remains to extend the result to all $\lambda$ on a common set of full
$\mathbb P$-measure. But, since the left-hand side in
\eqref{hyperplane_quenched} is nondecreasing in $\lambda$, as well as
the right-hand side is continuous in $\lambda$, such a set is given by
$\bigcap_{\lambda'\in[0,\infty)\cap\mathbb Q}\Omega_{\lambda'}$.

We step forward to the proof of Theorem~\ref{phase_transitions}.
Again, we restrict to the quenched setting. The annealed part of the theorem
then follows by a simple change of notations. 

We first want to establish \eqref{supremum_qu}. By Theorem~\ref{LDP}, it
suffices to show
\begin{equation}\label{supremum_qu*}
\sup_{x\in \mathbb R^d}\left(h\cdot x-I(x)\right)=
\left\{\begin{array}{cl}
0& \text{if $\alpha^*_0(h)\leq 1$}\,,\\
\lambda^{\mathrm{qu}}_h& \text{if $\alpha^*_0(h)> 1$}\,,
\end{array}\right.
\end{equation}
where $\lambda^{\mathrm{qu}}_h>0$  is the unique number with
$\alpha^{*}_{\lambda^{\mathrm{qu}}_h}(h)=1$.
Existence and uniqueness of
$\lambda^{\mathrm{qu}}_h$, as well as the property
$\lambda^{\mathrm{qu}}_h>0$, thereby follow from 
Lemma~\ref{properties_dual_norm}. 

In the case $\alpha^*_0(h)\leq 1$, the lower estimate for the
supremum is obvious. Assume now $\alpha^*_0(h)>1$.
From Theorem~\ref{shape_theorem}, we know that
 $\alpha_\lambda(x)$, to  fixed $x\in\mathbb R^d\setminus\{0\}$, is concave 
and strictly increasing in $\lambda\in[0,\infty)$. 
Therefore, the right
derivative 
\begin{align*}
\dot\alpha^+_{\lambda^{\mathrm{qu}}_h}(x)
=\lim_{\varepsilon\downarrow 0}\,
\frac{\alpha_{\lambda^{\mathrm{qu}}_h+\varepsilon}(x)
-\alpha_{\lambda^{\mathrm{qu}}_h}(x)}{\varepsilon}
\end{align*} 
is  well-defined and strictly
positive. From 
$\alpha_{\lambda^{\mathrm{qu}}_h}(x)$, it furthermore inherits
the scalar linearity in $x\in\mathbb R^d$. 
As a consequence, there exist $e\in\mathcal S^{d-1}$ and thus
$y=e/\dot\alpha^+_{\lambda^{\mathrm{qu}}_h}(e)$ with
$\dot\alpha^+_{\lambda^{\mathrm{qu}}_h}(y)=1$ and
$$
1=\alpha^*_{\lambda^{\mathrm{qu}}_h}(h)=\sup_{x\in\mathcal S^{d-1}}
\bigg(\frac{h\cdot x}{\alpha_{\lambda^{\mathrm{qu}}_h}(x)}\bigg)
=\frac{h\cdot e}{\alpha_{\lambda^{\mathrm{qu}}_h}(e)}
=\frac{h\cdot y}{\alpha_{\lambda^{\mathrm{qu}}_h}(y)}\,.
$$
By the first condition on $y$, the map 
$\lambda\mapsto\alpha_\lambda(y)-\lambda$  is nondecreasing for
$\lambda\leq\lambda^{\mathrm{qu}}_h$ and nonincreasing for
$\lambda>\lambda^{\mathrm{qu}}_h$. We therefore have
$I(y)=
\alpha_{\lambda^{\mathrm{qu}}_h}(y)-\lambda^{\mathrm{qu}}_h$. 
The second condition on $y$ now implies
\begin{equation*}
\sup_{x\in\mathbb R^d}\left(h\cdot x-I(x)\right)
\geq 
h\cdot y-I(y)
= \lambda^{\mathrm{qu}}_h\,.
\end{equation*}

For the reversed estimate,
we additionally set
$\lambda^{\mathrm{qu}}_h= 0$  when  $\alpha^*_0(h)\leq 1$.
We can assume
$h\neq 0$. The definition of $\alpha^*_\lambda$ then yields
\begin{align*}
\sup_{y:\, h\cdot y\geq 0}\left(h\cdot y-\alpha_{\lambda}(y)\right)
&\leq \sup_{y\in\mathbb R^d}
\left(h\cdot y-\frac{h\cdot y}{\alpha^{*}_{\lambda}(h)}\right) =
\left\{\begin{array}{cl}
0& \text{if $\lambda=\lambda^{\mathrm{qu}}_h$}\,,\\
\infty& \text{if $\lambda\neq\lambda^{\mathrm{qu}}_h$}\,,
\end{array}\right.
\intertext{which leads to}
\sup_{x\in\mathbb R^d}\left(h\cdot x-I(x)\right)
&\leq \inf_{\lambda\geq 0\vphantom{R^R_0}}\left(
\sup_{\;y:\, h\cdot y\geq 0}
\left(h\cdot y-\alpha_{\lambda}(y)\right)+\lambda\right)
=\lambda^{\mathrm{qu}}_h\,.
\end{align*}
This completes the proof of \eqref{supremum_qu*}.

It remains to establish the limiting behavior of $S(n)/n$. To this end, observe
that the rate function
$I_h$ satisfies
$$I_h(x)\geq\alpha(x)-h\cdot x$$
for all $x\in\mathbb R^d$.
When $\alpha^*_0(h)<1$, it only vanishes at the origin, and the
sub-ballistic behavior  follows by the large deviation
estimates in the quenched part of Theorem~\ref{LDP}.

On the other hand, when $\alpha^*_0(h)>1$, the rate function $I_h$ only
vanishes on the set  $M_h$, which is compact by the continuity of $I$  on 
its effective domain $D_I$ (which itself is compact).  Observe
furthermore that
$M_h$ cannot  contain the origin since
$\lambda^{\mathrm{qu}}_h>0$.
The ballistic behavior now follows again by the large deviation estimates.
This completes the proof of Theorem~\ref{phase_transitions}.


\begin{thebibliography}{9}


\bibitem[Deuschel, Stroock(1989)]{DS}
J.D.~Deuschel, D.W~Stroock,
\textit{Large deviations.} 
Academic Press, San Diego 1989.

\bibitem[Flury(2006)]{Flu}
M.~Flury,
Coincidence of Lyapunov exponents for random walks in weak random potentials,
http://arxiv.org/abs/math.PR/0608357 (2006).

\bibitem[Grimmett(1999)]{Gri}
G.~Grimmett,
\textit{Percolation.}  
Springer-Verlag, Berlin 1999.

\bibitem[Sznitman(1995a)]{Sz95a}
A.S.~Sznitman, 
Quenched critical large deviations for brownian motion in a Poissonian potential. 
\textit{J.\ Funct.\ Anal.\ {\bf131}} (1995a) 45--77.

\bibitem[Sznitman(1995b)]{Sz95b}
A.S.~Sznitman, 
Annealed Lyapunov exponents and large deviations in
a Poissonian potential. II. 
\textit{Ann.\ Sci.\ Ecole Norm.\ Sup.\ (4) {\bf28}} (1995b) 371--390.

\bibitem[Sznitman(1998)]{Szn}
 A.S.~Sznitman,
\textit{Brownian motion, obstacles, and random media.}
Springer, Berlin 1998.

\bibitem[Woess(2000)]{Woe}
 W.~Woess,
\textit{Random walks on infinite graphs and groups.}
Cambridge University Press, Cambridge 2000.

\bibitem[Zerner(1996)]{Zer}
M.P.W.~Zerner,
Directional decay of the Green's function for a
random nonnegative potential on $\mathbb Z^d$,
\textit{Ann.\ Appl.\ Probab.\ {\bf 8}} (1996) 246--280.

\end{thebibliography}
\end{document}